\documentclass[review]{elsarticle}
\usepackage[top=32mm, bottom=32mm, left=25mm, right=25mm]{geometry}
\usepackage[colorlinks=true]{hyperref}
\usepackage{lineno,amssymb,amsmath,latexsym,mathrsfs}
\usepackage{graphicx}
\usepackage[all]{xy}
\usepackage{amssymb,amsbsy}
\usepackage[active]{srcltx} 
\usepackage{graphics}
\usepackage{eucal}
\usepackage{amsmath}
\usepackage{graphicx}
\usepackage{caption}
\usepackage{subcaption}
\usepackage{enumerate}
\journal{}
\begin{document}

\newtheorem{theorem}{Theorem}[section]
\newtheorem{algorithm}[theorem]{Algorithm}
\newtheorem{axiom}[theorem]{Axiom}
\newtheorem{case}[theorem]{Case}
\newtheorem{claim}[theorem]{Claim}
\newtheorem{conclusion}[theorem]{Conclusion}
\newtheorem{condition}[theorem]{Condition}
\newtheorem{conjecture}[theorem]{Conjecture}
\newtheorem{cor}[theorem]{Corollary}
\newtheorem{criterion}[theorem]{Criterion}
\newtheorem{definition}[theorem]{Definition}
\newtheorem{example}[theorem]{Example}
\newtheorem{exercise}[theorem]{Exercise}
\newtheorem{lemma}[theorem]{Lemma}
\newtheorem{notation}[theorem]{Notation}
\newtheorem{corollary}[theorem]{Corollary}
\newtheorem{proposition}[theorem]{Proposition}
\newtheorem{remark}[theorem]{Remark}
\newtheorem{solution}[theorem]{Solution}
\newtheorem{summary}[theorem]{Summary}
\newcommand{\pf}{\noindent \textbf{Proof.}\quad}
\newcommand{\epf}{\hspace{\stretch{1}}$\blacksquare$}

\begin{frontmatter}

\title{Some Properties of Recurrent Sets for Endomorphisms\\ of  Topological Groups}

\author[a1]{Seyyed Alireza Ahmadi}
\address[a1]{Department of Mathematics, University of Sistan and Baluchestan, Zahedan, Iran.}
\ead{sa.ahmadi@math.usb.ac.ir, sa.ahmdi@gmail.com}

\author[a1]{Javad Jamalzadeh}
\ead{jamalzadeh1980@math.usb.ac.ir}

\author[a2]{Xinxing Wu\corref{mycorrespondingauthor}}
\cortext[mycorrespondingauthor]{Corresponding author}
\address[a2]{School of Sciences, Southwest Petroleum University, Chengdu, Sichuan 610500, People's
Republic of China}
\ead{wuxinxing5201314@163.com}

\begin{abstract}
This paper studies topological definitions of chain recurrence and shadowing for continuous endomorphisms of topological
groups generalizing the relevant concepts for metric spaces. It is proved that in this case the sets of chain recurrent
points and chain transitive component of the identity are topological subgroups. Furthermore, it is obtained that some
dynamical properties induced by the original system on quotient spaces. These results link an algebraic property to a dynamical property.
\end{abstract}
\end{frontmatter}


\section{Introduction}
One of the main problems in discrete and continuous dynamical systems is the description of orbit structure for a system from a topological
point of view. A discrete dynamical system usually consists of a compact metric space $X$ and a continuous function $f$ from $X$ to itself.
A number of properties of interest in such systems are defined in purely topological terms, for example recurrence, non-wandering points
(see below for details). Recently, Good and Macias \cite{good2018topological} defined other properties for dynamical systems in purely
topological terms, for example sensitive dependence on initial conditions, chain transitivity and recurrence, shadowing, and positive expansiveness.
In the presence of compactness, existence of an unique uniformity,  allows ue to mimic existing metric proofs. The uniform
approach has been studied in a number of cases: Hood~\cite{MR0353282} defined topological entropy for uniform spaces; Morales and
Sirvent~\cite{MR3458385} considered positively expansive measures for measurable functions on uniform spaces, extending results from the
literature; Devaney chaos for uniform spaces is considered in \cite{MR3085695}; Auslander et al.~\cite{MR3223369} generalized many
known results about equicontinuity to the uniform spaces; Das et al. \cite{Das2013149} generalized spectral decomposition theorem to the uniform spaces;
{\color{blue}We \cite{AWFML2018} generalized  concepts of entropy points, expansivity and shadowing
property for dynamical systems on uniform spaces and we obtained a relation between
topological shadowing property and positive uniform entropy;}; Wu et al. \cite{WMZL2018} obtained that every point transitive dynamical system defined on a Hausdorff uniform space is either almost
equicontinuous or sensitive.

Motivated by these ideas we show that if  the underlying set of a dynamical system is an abelian topological group, then surprisingly dynamical objects of a dynamical system exhibit some algebraic properties. Recurrence behavior  is one of the most important concepts in topological dynamics \cite{MR2330357,SILVESTROV2002117,Wu2018}. We are going to investigate the properties of recurrent sets of a continuous endomorphism of a topological group as a discrete dynamical system.
A topological group is a set $G$ on which two structures are given, a group structure and a topology, such that the group operations are continuous. Specifically, the mapping $(x,y)\mapsto xy^{-1}$ from the direct product $G\times G$ into $G$ must be continuous. A subgroup of a topological group is a topological group in the induced topology. Every topological group is a uniform space in a natural way. Specifically, a  uniform group structure on a topological group is defined by the collection of sets
$$
\left\{(x,y)~|~xy^{-1}\in E\right\};\quad E\in\mathfrak{B}_e,
$$
where $\mathfrak{B}_e$ is a system of symmetric neighborhoods of the identity $e$ in $G$. Make a standardized assumption that all topological groups are abelian and compact, and $f$ is a continuous endomorphism on $G$, although some of the results apply to more general settings. A fixed point of dynamical system $f$, exhibits the simplest type of recurrence. We denote by $Fix(f)$ the set of all fixed points of $f$.

A point carried back to itself by a dynamical system $f$ exhibits the next most elementary type of recurrence. For some $m\in \mathbb{N}$, a point $x\in G$ is
called {\it $m$-periodic} if $f^m(x)=x$. We denote by $Per_m(f)$ the set of all $m$-periodic points of $f$ and we set $Per(f)=\bigcup_{m=1}^{\infty}Per_m(f)$.
A point $x\in G$ is {\it non-wandering} if for each neighborhood $U$ of $x$, there exists $n\in\mathbb{N}$ such that $U\cap f^n(U)\neq \O$. We denote by
$\Omega(f)$ the set of all non-wandering points of $f$.

For $D\in\mathcal{B}_e$, a {\it $D$-pseudo-orbit} or {\it $D$-chain} of $f$ is a sequence $\{x_n\}_{n\in \mathbb{N}}$ such that $f(x_n)x_{n+1}^{-1}\in D$
for $n\in\mathbb{N}$. We use the symbol $\mathcal{O}_{E}(f,x,y)$ for the set of $E$-chains $\{x_0,x_1,\dots, x_n\}$ of $f$ with $x_0=x$ and $x_n=y$. For
$x, y\in G$, we write  $x\stackrel{E}{\rightsquigarrow} y$ if $\mathcal{O}_{E}(f,x,y)\neq\O$ and we write $x\rightsquigarrow y$ if
$\mathcal{O}_{E}(f,x,y)\neq\O$ for each $E\in\mathfrak{B}_e$. We write $x\leftrightsquigarrow y$ if $x\rightsquigarrow y$ and
$y\rightsquigarrow x$. The set $\{x\in G~| ~x\leftrightsquigarrow x\}$ is called the \textit{chain recurrent set} of $f$ and denoted by $\mathrm{CR}(f)$.
Denote by $\mathrm{CC}(f)$ the \textit{chain component} of $f$ containing the identity $e$, i.e., $\mathrm{CC}(f)=\{x\in G~|~ x\leftrightsquigarrow e\}$ \cite{GARAY1989372}. Clearly,
$$
Fix(f)\subseteq Per_m(f)\subseteq Per(f)\subseteq \Omega(f)\subseteq \mathrm{CR}(f).
$$

\section{Recurrent Subgroups}
This section is devoted to the algebraic properties of recurrent sets. Our following results show that when the underlying set of a dynamical system is a topological group, then most of the well-known recurrent sets are also topological subgroups of the underlying topological group.
We seek a definition of recurrence so that the set $\mathfrak{R}(f)$ of recurrent points
with respect to an endomorphism $f$ has the following desirable properties:
\begin{enumerate}[(R1)]
\item $\mathfrak{R}(f)$ is a subgroup;
\item
 The set $\mathfrak{R}(f)$ is forward invariant with respect to $f$, i.e., $f(\mathfrak{R}(f))\subseteq\mathfrak{R}(f)$.
\item $\mathfrak{R}(f)$ is closed;
\item
$\mathfrak{R}(f)$ is invariant under topological conjugacy, i.e., if $f:G\rightarrow G$ and $g:H\rightarrow H$ are two continuous endomorphisms on topological groups and $\phi:G\rightarrow H$ is a continuous automorphism with continuous inverse such that $\phi\circ f=g\circ \phi$, then $\mathfrak{R}(g)=\phi(\mathfrak{R}(f))$;
\item
$\mathfrak{R}(f)$ is invariant under canonical mapping, i.e., if $\tilde{f}:G/H \rightarrow G/H$ is the canonical mapping induced by $f$, then $\mathfrak{R}(\tilde{f})=\{\mathfrak{R}(f)\}$.
\end{enumerate}
\begin{definition}
We say that a subset $\mathfrak{R}(f)$ of $G$ with respect to an endomorphism $f$ is a \textit{recurrent subgroup} if it satisfies properties  (R1)--(R4).
\end{definition}

Closely related to fixed points are the eventually fixed points, which are the points that reach a fixed point after finitely many iterations.
More explicitly, a point $x$ is said to be an {\it eventually fixed point} of a map $f$ if there exists some $k\in \mathbb{N}$ such that
$f^k(x)\in Fix(f)$. A point $x$ is said to be an {\it eventually $m$-periodic point} if $f^{k}(x) \in Per_m(f)$ for some $k\in \mathbb{N}$.
Denote by $EFix(f)$ and $EPer_m(f)$, the set of all eventually fixed point and the set of eventually $m$-periodic points of $f$, respectively.
Also we set $EPer(f)=\bigcup_{m=1}^{\infty}EPer_m(f)$.

\begin{proposition}\label{p1}
Let $f$ be a continuous endomorphism of a topological group $G$. Then
\begin{enumerate}[(i)]
\item
$Fix(f)$ is a recurrent subgroup.
\item
$Per_m(f)$ is a recurrent subgroup.
\item
$\overline{Per(f)}$ is a recurrent subgroup.
\item
$\overline{EFix(f)}$ is a recurrent subgroup.
\item
$\overline{EPer_m(f)}$ is a recurrent subgroup.
\item
$\overline{EPer(f)}$ is a recurrent subgroup.
\end{enumerate}
\end{proposition}

\pf
(i) For any $x,y\in Fix(f)$, it can be verified that $f(xy^{-1})=f(x)f(y^{-1})=xy^{-1}$, implying that $Fix(f)$ is a subgroup.
Clearly, $Fix(f)$ is closed, $f$-invariant and invariant under algebraic topological conjugacy.

(ii) The proof is similar to part (i).

(iii) For any $x, y\in Per(f)$, there exist $m, n\in \mathbb{N}$ such that $f^m(x)=x$ and $f^n(y)=y$, implying that $f^{mn}(xy^{-1})=(f^m)^n(x)(f^n)^m(y^{-1})=xy^{-1}$. Then, $xy^{-1}\in Per(f)$. This implies that $Per(f)$ and so $\overline{Per(f)}$ is a subgroup.
Clearly, $Per(f)$ and so by continuity $\overline{Per(f)}$ is $f$-invariant. Let $f:G\rightarrow G$ and $g:H\rightarrow H$ be two continuous endomorphisms
and let $\phi:G\rightarrow H$ be a continuous automorphism with continuous inverse such that $\phi\circ f=g\circ \phi$. Clearly, $\phi(Per(f))\subset Per(g)$.
This implies that $\phi(\overline{Per(f)})\subset\overline{\phi(Per(f))}\subset \overline{Per(g)}$.
For the reverse inclusion, let $x\in \overline{Per(g)}$. Then there exists a net $x_{\lambda}$ in $Per(g)=\phi(Per(f))$ with $x_{\lambda}\rightarrow x$.
For each $\lambda$, there exists a point $z_{\lambda}\in Per(f)$ such that $x_{\lambda}=\phi(z_{\lambda})$. This, together with $x_{\lambda}\rightarrow x$,
implies that $z_{\lambda}\rightarrow\phi^{-1}(x)$. Therefore, $\phi^{-1}(x)\in\overline{Per(f)}$. This implies that $\phi(\overline{Per(f)})= \overline{Per(g)}$.

(iv)
For any $x,y\in EFix(f)$, there exist $m, n\in \mathbb{N}$ such that $f^m(x),f^n(y)\in Fix(f)$, implying that
$f^{mn}(xy^{-1})=f^m(x)f^n(y^{-1})\in Fix(f)$. Then, $xy^{-1}\in EFix(f)$. This implies that $EFix(f)$ and so $\overline{EFix(f)}$ is a subgroup.
Clearly, $EFix(f)$ and so by continuity $\overline{EFix(f)}$ is $f$-invariant. Let $f:G\rightarrow G$ and $g:H\rightarrow H$ be two continuous
endomorphisms and let $\phi:G\rightarrow H$ be a continuous automorphism with continuous inverse such that $\phi\circ f=g\circ \phi$. If
$\phi(x)\in \phi(EFix(f))$, then there exists positive integer $m$ such that $f^m(x)\in Fix(f)$. Thus, $\phi(f^m(x))\in\phi(Fix(f))=Fix(g)$,
implying that $g^m(\phi(x)\in Fix(g)$. Therefore, $\phi(EFix(f))\subset EFix(g)$. This implies that
$\phi(\overline{EFix(f)})\subset\overline{\phi(EFix(f))}\subset \overline{EFix(g)}$. For the reverse inclusion, let $x\in \overline{EFix(g)}$.
Then there exists a net $x_{\lambda}$ in $EFix(g)=\phi(EFix(f))$ with $x_{\lambda}\rightarrow x$. For each $\lambda$, there exists a point
$z_{\lambda}\in EFix(f)$ such that $x_{\lambda}=\phi(z_{\lambda})$. This, together with $x_{\lambda}\rightarrow x$, implies that $z_{\lambda}\rightarrow\phi^{-1}(x)$.
Therefore, $\phi^{-1}(x)\in\overline{EFix(f)}$. This implies that $\phi(\overline{EFix(f)})= \overline{EFix(g)}$.

(v)
For any $x, y\in Eper_m(f)$, there exist $k, l\in \mathbb{N}$ such that $f^{m+k}(x)=f^k(x)$ and $f^{m+l}(y)=f^l(y)$, implying that $f^{m+k+l}(x)=f^{k+l}(x)$ and $f^{m+k+l}(y)=f^{k+l}(y)$. Then, $f^{m+k+l}(xy^{-1})=f^{k+l}(xy^{-1})$. This implies that $EPer_m(f)$ and so $\overline{EPer_m(f)}$ is a subgroup.
Similarly, one can prove the properties (R3) and (R4) by adapting the proof of part (iii).

(vi)
For any $x, y\in Eper(f)$, there exist $k, l, m\in \mathbb{N}$ and $n\in \mathbb{N}$ such that
$f^{m+k}(x)=f^k(x)$ and $f^{n+l}(y)=f^l(y)$, implying that $f^{m+k+l}(x)=f^{k+l}(x)$ and $f^{n+k+l}(y)=f^{k+l}(y)$.
Applying induction yields that
$$
f^{2m+k+l}(x)=f^m(f^{m+k+l}(x))=f^{m+k+l}(x)=f^{k+l}(x),
$$
$$
 \vdots
$$
$$
f^{mn+k+l}(x)=f^{m(n-1)+k+l}(x)=f^{k+l}(x),
$$
 and
$$
f^{2n+k+l}(y)=f^n(f^{n+k+l}(y))=f^{n+k+l}(y)=f^{k+l}(y),
$$
$$
 \vdots
$$
$$
f^{mn+k+l}(y)=f^{n(m-1)+k+l}(y)=f^{k+l}(y).
$$

Thus, $f^{mn+k+l}(xy^{-1})=f^{k+l}(xy^{-1})$. This implies that $EPer(f)$ and so $\overline{EPer(f)}$ is a subgroup.
Similarly, one can prove the properties(R3) and (R4) by adapting the proof of part (iii).
\epf

The following proposition shows that a chain recurrent set, which includes all the types of
returning trajectories: periodic, eventually-periodic, non-wandering and so on, is a subgroup of $G$.

\begin{proposition}
Let $f$ be a continuous endomorphism of topological group $G$. Then $\mathrm{CR}(f)$ is a recurrent subgroup.
\end{proposition}
\pf
For any $x\in \mathrm{CR}(f)$ and any $E\in\mathfrak{B}_e$, there exists an $E$-chain $\{x_1,x_2,\dots,x_n\}$ with $x_1=x_n=x$ such that
$f(x_{n})x_{n+1}^{-1}\in E$ for all $1\leq i\leq n-1$, implying that $f(x_i^{-1})x_{i+1}\in E^{-1}=E$. Then, $x^{-1}\in \mathrm{CR}(f)$.

For any $x,y\in \mathrm{CR}(f)$, we show that $xy\in \mathrm{CR}(f)$. In fact, for any $E\in\mathfrak{B}_e$, take some $W\in\mathfrak{B}_e$
such that $W^2\subset E$. Since $x,y\in \mathrm{CR}(f)$, there exist $E$-chains $\{x_1,x_2,\dots,x_m\}$ and $\{y_1,y_2,\dots,y_n\}$ with $x_1=x_m=x$ and $y_1=y_n=y$ such that $f(x_i)x_{i+1}^{-1}\in W$ for $1\leq i\leq m-1$ and $f(y_j)y_{j+1}\in W$ for $1\leq j\leq n-1$. Choose two extended sequences $\{x_i\}_{i=1}^{mn}$ and $\{y_j\}_{j=1}^{mn}$ as following:
\begin{eqnarray*}
&&x_{i+kn}=x_i \text{ for } 1\leq i\leq m, \ 0\leq k\leq n-1;\\
&&y_{j+kn}=y_j \text{ for } 1\leq j\leq n, \ 0\leq k\leq m-1.
\end{eqnarray*}
Clearly, $x_1y_1=x_{mn}y_{mn}=xy$ and
$$
f(x_iy_i)(x_{i+1}y_{i+1})^{-1}=f(x_i)x_{i+1}^{-1}f(y_i)y_{i+1}^{-1}\in W^2\subset E.
$$
Therefore, $xy\in \mathrm{CR}(f)$, implying that $\mathrm{CR}(f)$ is a subgroup.

Next assume that $E\in\mathfrak{B}_e$ and choose $\hat{E}\in\mathfrak{B}_e$ such that $\hat{E}^2\subset E$. By uniform continuity there exists
$D\in\mathfrak{B}_e$ such that $xy^{-1}\in D$ implies $f(x)f(y)^{-1}\in\hat{E}$. Choose $\hat{D}\in\mathfrak{B}_e$ with $\hat{D}^2\subset D$.
Assume that $x^{(\lambda)}$ be a net in $\mathrm{CR}(f)$ such that $x^{(\lambda)}\rightarrow x$. Then for some $\lambda$, $x^{(\lambda)}x^{-1}\in D$.
Since $x^{(\lambda)}\in \mathrm{CR}(f)$, there exists a $\hat{D}$-pseudo-orbit $\{x_0,x_1,\dots,x_n\}$ with $x_0=x_n=x^{(\lambda)}$. Thus
$\{x,x_1,\dots,x_{n-1},x\}$ is an $E$-pseudo-orbit and hence $x\in \mathrm{CR}(f)$. Therefore, $\mathrm{CR}(f)$ is closed.

Fix any $E\in\mathfrak{B}_e$. Then, $\{x,f(x)\}$ is an $E$-pseudo-orbit from $x$ to $f(x)$. Choose $U,V\in\mathfrak{B}_e$ such that $U^2\subset E$ and
$V\subset U\cap f^{-1}(U)$. Since $x\in \mathrm{CR}(f)$, there exists a $V$-pseudo-orbit $\{x_0=x,x_1,\dots,x_{n-1},x_n=x\}$ from $x$ to itself. Then $\{f(x),x_2,x_3,\dots,x_n=x\}$ is an $E$-pseudo-orbit from $f(x)$ to $x$ and hence $\{f(x),x_2,x_3,\dots,x_{n-1},x,f(x)\}$ is an $E$-pseudo-orbit
from $f(x)$ to itself. Therefore, $f(x)\in \mathrm{CR}(f)$, implying that $f(\mathrm{CR}(f))\subset \mathrm{CR}(f)$.
\epf

The relation `$\leftrightsquigarrow$' is an equivalence relation on $\mathrm{CR}(f)$.
The equivalence classes of this relation are called {\it chain components}. These are compact invariant sets and cannot
be decomposed into two disjoint nonempty compact invariant sets, hence serve as building
blocks of the dynamics. The topology of chain recurrent set and chain components
have been always in particular interest~\cite{ barros2010attractors, bonatti2014tame, MR2379487, Shekutkovski2015}.
\begin{proposition}
Let $f$ be a continuous endomorphism of topological group $G$. Then $\mathrm{CC}(f)$ is a recurrent subgroup.
\end{proposition}
\pf
Suppose that $x,y\in \mathrm{CC}(f)$ and $E\in\mathfrak{B}_e$. Choose $W\in\mathfrak{B}_e$ such that $W^2\subset E$. Then,
$x\stackrel{W}{\rightsquigarrow} e$ and $y\stackrel{W}{\rightsquigarrow} e$, implying that there exist $W$-chains $\{x=x_0,x_1,\dots,x_m=e\}$
and $\{y=y_0,y_1,\dots,y_n=e\}$. Without loss of generality, assume that $m\leq n$. Clearly,
$\{x_0^{-1}y_0,x_1^{-1}y_1,\dots x_m^{-1}y_m,y_{m+1},\dots,y_n=e\}$ is an $E$-chain from $x^{-1}y$ to $e$. It follows from
$x,y\in \mathrm{CC}(f)$ that $e\stackrel{W}{\rightsquigarrow} x$ and $e\stackrel{W}{\rightsquigarrow} y$. Then, there exist $W$-chains
$\{e=\hat{x}_0,\hat{x}_1,\dots,\hat{x}_p=x\}$ and $\{e=\hat{y}_0,\hat{y}_1,\dots,\hat{y}_q=y\}$. Without loss of generality, assume that $p\leq q$.
Then, the sequence $\{\hat{y}_0,\hat{y}_1,\dots,\hat{y}_{q-p-1},\hat{y}_{q-p}\hat{x}_0^{-1},\hat{y}_{q-p+1}\hat{x}_1^{-1},\dots,\hat{y}_q \hat{x}_p^{-1}\}$
is an $E$-chain from $e$ to $x^{-1}y$. Therefore, $x^{-1}y\in \mathrm{CC}(f)$, implying that $\mathrm{CC}(f)$ is a subgroup.

Let $E\in\mathfrak{B}_e$ and choose $W\in\mathfrak{B}_e$ such that $W^2\subset E$. By uniform continuity there exists $W\supset D\in\mathfrak{B}_e$
such that $xy^{-1}\in D$ implies $f(x)f(y^{-1})\in W$. Let $z\in\overline{\mathrm{CC}(f)}$. Then, there exists $x\in \mathrm{CC}(f)$ such that $xz^{-1}\in D$.
Clearly, $x\stackrel{D}{\rightsquigarrow} e$ and $e\stackrel{D}{\rightsquigarrow} x$. This implies that there exist $D$-chains $\{x=x_0,x_1,\dots,x_n=e\}$
and $\{e=x'_0,x'_2,\dots,x'_m=x\}$. Clearly, $\{z,x_1,x_2,\dots,e\}$ and $\{e=x'_0,x'_1,\dots,x'_{m-1},z\}$ are $E$-chains. Then,
$z\stackrel{E}{\rightsquigarrow} e$ and $e\stackrel{E}{\rightsquigarrow} z$. Therefore, $z\in \mathrm{CC}(f)$, implying that $\mathrm{CC}(f)$ is closed.

Suppose that $x\in \mathrm{CC}(f)$ and $E\in\mathfrak{B}_e$. Choose $U,V\in\mathfrak{B}_e$ such that $U^2\subset E$ and $V\subset U\cap f^{-1}(U)$.
It follows from $x\in \mathrm{CC}(f)$ that there exists a $V$-pseudo-orbit $\{x_0=x,x_1,\dots,x_{n-1},x_n=e\}$ from $x$ to $e$. Then,
$\{f(x),x_2,x_3,\dots,x_n=e\}$ is an $E$-pseudo-orbit from $f(x)$ to $e$. Again from $x\in \mathrm{CC}(f)$, it follows that there exists a
$V$-pseudo-orbit $\{x_0=e,x_1,\dots,x_{n-1},x_n=x\}$ from $e$ to $x$. Then, $\{x_0,x_1,x_2,x_3,\dots,x_n=x,f(x)\}$ is an $E$-pseudo-orbit
from $e$ to $f(x)$. Thus, $f(x)\in \mathrm{CC}(f)$, implying that $f(\mathrm{CC}(f))\subset \mathrm{CC}(f)$.
\epf
\begin{remark}
It is not difficult to check that if $f$ is an automorphism, then $f(\mathrm{CC}(f))=\mathrm{CC}(f)$.
\end{remark}
\section{Dynamics Induced on Quotient Spaces by endomorphisms}

Suppose that $G$ is a topological group with identity $e$, and $H$ is a closed subgroup of
$G$. Denote by $G/H$ the set of all left cosets $aH$ of $H$ in $G$, and endow it with the
quotient topology with respect to the canonical mapping $\pi:G\rightarrow G/H$ defined by
$\pi(x)=xH$ for any $x\in G$. Then, the family $\{\pi(xE)~|~x\in G,~E\in\mathfrak{B}_e\}$
is a local base of the space $G/H$ at the point $xH\in G/H$, the mapping $\pi$ is open, and
$G/H$ is a homogeneous $T_1$-space. An endomorphism $f:G\rightarrow G$ induced a map
$\tilde{f}:G/H\rightarrow G/H$ such that the following diagram is commuted:
\begin{displaymath}
\xymatrix{
G \ar[r]^{f} \ar[d]_{\pi} &
G \ar[d]^{\pi} \\
G/H \ar[r]_{\tilde{f}} & G/H }
\end{displaymath}
This mapping is called {\it canonical map} \cite{bagley_peyrovian_1986}. In last section we introduce several recurrent
subgroups $\mathfrak{R}(f)$ for a dynamical system $f$, which leads us to investigate the dynamic of induced mapping
$\tilde{f}:G/\mathfrak{R}(f)\rightarrow G/\mathfrak{R}(f)$. We are interested in cases that $\mathfrak{R}(f)$ is invariant
under canonical mapping, i.e., $\mathfrak{R}(\tilde{f})=\{\mathfrak{R}(f)\}$.

\begin{proposition}
Let $f:G\rightarrow G$ be a continuous endomorphism and $\tilde{f}:G/Fix(f)\rightarrow G/Fix(f)$ be the
map induced by $f$. Then, $Fix(\tilde{f})=\{Fix(f)\}$.
\end{proposition}
\pf
Suppose that $H=Fix(f)$ and $\pi(x)=xH\in Fix(\tilde{f})$. Then, $\tilde{f}(xH)=xH$. Thus, $f(x)x^{-1}\in H=Fix(f)$.
Therefore, $f^n(f(x)x^{-1})=f(x)x^{-1}$ for any $n\in\mathbb{N}$. Applying induction implies that $f^{n+1}(x)=(f(x))^{n+1}x^{-1}$.
Hence,
\begin{align*}
 &\left\lbrace
  \begin{array}{c l}
  &f^2(x)=f(x)^2x^{-1}     \vspace{2pt}\Rightarrow f^3(x)=f^2(x)f(x)\Rightarrow f^4(x)=f^3(x)f^2(x),\\
  &f^3(x)=f(x)^3x^{-1},\\
  \end{array}
\right.\\
&\left\lbrace
  \begin{array}{c l}
  &f^3(x)=f(x)^3x^{-1}     \vspace{2pt}\Rightarrow f^4(x)=f^3(x)f(x),\\
  &f^4(x)=f(x)^4x^{-1}.\\
  \end{array}
\right.
\end{align*}
Thus $f(x)=f^2(x)=(f(x))^2x^{-1}$, implying that $f(x)=x$. Therefore, $xH=H$.
\epf
\begin{proposition}
Let $f:G\rightarrow G$ be a continuous endomorphism and $\tilde{f}:G/Per_m(f)\rightarrow G/Per_m(f)$ be the
map induced by $f$. Then for any $m\in \mathbb{N}$, $Per_m(\tilde{f})=\{Per_m(f)\}$.
\end{proposition}
\pf
 Let $H=Per_m(f)$ and $\pi(x)=xH\in Per_m(\tilde{f})$. Then, $\tilde{f}^m(xH)=xH$. Thus, $f^m(x)x^{-1}\in H=Per_m(f)$.
 Therefore, $f^m(f^m(x)x^{-1})=f^{m}(x)x^{-1}$. Applying induction implies that $f^{im}(x)=(f^{m}(x))^ix^{-1}$ for any $i\geq 2$.
 Hence,
\begin{align*}
 &\left\lbrace
  \begin{array}{c l}
  &f^{2m}(x)=f^m(x)^{2}x^{-1}     \vspace{2pt}\Rightarrow f^{3m}(x)=f^{2m}(x)f^m(x)\Rightarrow f^{4m}(x)=f^{3m}(x)f^{2m}(x),\\
  &f^{3m}(x)=f^m(x)^{3}x^{-1},\\
  \end{array}
\right.\\
&\left\lbrace
  \begin{array}{c l}
  &f^{3m}(x)=f^m(x)^{3}x^{-1}     \vspace{2pt}\Rightarrow f^{4m}(x)=f^{3m}(x)f^m(x),\\
  &f^{4m}(x)=f^m(x)^{4}x^{-1}.\\
  \end{array}
\right.
\end{align*}
Thus, $f^m(x)=f^{2m}(x)=(f^m(x))^2x^{-1}$, implying that $f^m(x)=x$.
\epf
\begin{proposition}
Let $f:G\rightarrow G$ be a continuous automorphism and $\tilde{f}:G/\mathrm{CC}(f)\rightarrow G/\mathrm{CC}(f)$
be the map induced by $f$. Then, $\mathrm{CC}(\tilde{f})=\{\mathrm{CC}(f)\}$.
\end{proposition}
\pf
Suppose that $H=\mathrm{CC}(f)$. Let $xH\in\mathrm{ CC}(\tilde{f})$ and $E\in\mathfrak{B}_e$. Then, $xH\stackrel{\pi(E)}{\rightsquigarrow} H$,
implying that there exist points $x_0,x_1,\dots,x_n\in G$ such that $x_0=x$,  $x_n\in H$ and $\tilde{f}(\pi(x_i))\pi(x_{i+1})^{-1}\in \pi(H)$
for all $0\leq i\leq n-1$. Thus, for any $1\leq i\leq n$, there exist $e_{i}\in E$ and $h_i\in H$ such that $f(x_i)x_{i+1}^{-1}=e_{i+1}h_{i+1}$.
Choose $y_i=x_ih'_i$ with
$$
h'_0=e,\quad  h'_{i+1}=f(h'_i)h_{i+1}.
$$
Then, $y_0=x$, $y_n=h'_n$, and
$$
f(y_i)y_{i+1}^{-1}=f(x_i)x_{i+1}^{-1}f(h'_i)(h'_{i+1})^{-1}=f(x_i)x_{i+1}^{-1}h_{i+1}^{-1}=e_{i+1}\in E,
$$
implying that $x\stackrel{E}{\rightsquigarrow} h'_n$. From $h_n'\in H=\mathrm{CC}(f)$, it follows that $h'_n\stackrel{E}{\rightsquigarrow} e$.
Then, $x\stackrel{E}{\rightsquigarrow} e$. This implies that $x{\rightsquigarrow} e$ due to the arbitrariness of $E$.

From $xH\in \mathrm{CC}(\tilde{f})$ and $H\stackrel{\pi(H)}{\rightsquigarrow} xH$, it follows that there exist points $x_0,x_1,\dots,x_n\in G$
such that $x_0\in H$, $x_n=x$ and $\tilde{f}(\pi(x_i))\pi(x_{i+1})^{-1}\in \pi(H)$ for all $0\leq i\leq n-1$. This implies that there exist
$e_{i}\in E$ and $h_i\in H$ such that $f(x_i)x_{i+1}^{-1}=e_{i+1}h_{i+1}$. For any $0\leq i\leq n$, choose $y_i=x_ih'_i$ with
$$
h'_n=e,\quad  h'_{i-1}=f^{-1}(h'_i)h_{i+1}^{-1}.
$$
Then, $y_0=h'_0$, $y_n=x$, and
$$f(y_i)y_{i+1}^{-1}=f(x_i)x_{i+1}^{-1}f(h'_i)(h'_{i+1})^{-1}=f(x_i)x_{i+1}^{-1}h_{i+1}^{-1}=e_{i+1}\in E,
$$
implying that $h_0'\rightsquigarrow x$. From $h_0'\in H=\mathrm{CC}(f)$, it follows that $e\stackrel{E}{\rightsquigarrow} h'_0$. Therefore, $e\stackrel{E}{\rightsquigarrow} x$, implying that $e{\rightsquigarrow} x$ due to the arbitrariness of $E$. Hence, $x\in \mathrm{CC}(f)$.
\epf

The shadowing property provides tools for fitting real trajectories nearby to approximate trajectories
\cite{Carvalho2014801}. The following definition generalizes the relevant concept for metric spaces to topological groups.
\begin{definition}
We say that a $D$-pseudo orbit of $f$ is {\it $E$-shadowed} by a point $x$ in $G$ if $f^n(x)x_n^{-1}\in E$ for any $n\in\mathbb{N}$.
A continuous endomorphism $f:G\rightarrow G$ has the {\it shadowing property} if for any $E\in\mathfrak{B}_e$, there exists
some $D\in\mathfrak{B}_e$ such that every $D$-pseudo orbit of $f$ can be $E$-shadowed by some point in $G$.
\end{definition}

\begin{lemma}\label{lem1}
Let $f:G\rightarrow G$ be a continuous endomorphism with the shadowing property. If $H$ is an $f$-invariant subgroup of $G$,
then for any $E\in\mathfrak{B}_e$, there exists $D\in\mathfrak{B}_e$ such that every $DH$-pseudo-orbit can be $EH$-shadowed
by some point in $G$.
\end{lemma}
\pf
Fix any $E\in\mathfrak{B}_e$. The shadowing property implies that there exists $D\in\mathfrak{B}_e$ such that every $D$-pseudo-orbit
can be $E$-shadowed by some point in $G$. Given any fixed $DH$-pseudo-orbit $\{y_n\}$, then $f(y_n)y_{n+1}^{-1}\in DH$ for all
$n\in \mathbb{N}$. This implies that for any $n\in\mathbb{N}$, there exists $d_n\in D$ and $h_n\in H$ such that $f(y_n)y_{n+1}^{-1}h_n^{-1}=d_n$.
Choose the sequence $h_n'$with $h_{n+1}'=f(h_n')h_n$ and take $x_n=y_nh_n'$. Then,
\begin{eqnarray*}
f(x_n)x_{n+1}&=&f(y_nh_n')(y_{n+1}h_{n+1}')^{-1}\\
&=&f(y_n)y_{n+1}^{-1}f(h_n')(h'_{n+1})^{-1}\\
&=&f(y_n)y_{n+1}h_n^{-1}=d_n\in D,
\end{eqnarray*}
implying that $\{x_n\}$ is a $D$-pseudo-orbit of $f$. Thus, there exists $y\in G$ such that $f^n(y)x_n^{-1}\in E$ for $n=0,1,2,\dots$.
Therefore, $f^n(y)(y_nh'_{n})^{-1}\in E$. This implies that $f^n(y)y_n^{-1}\in EH$.
\epf
\begin{theorem}
Let $f:G\rightarrow G$ be a continuous endomorphism  with the shadowing property and $H$ be a an
$f$-invariant subgroup of $G$. Then, the canonical mapping $\tilde{f}:G/H\rightarrow G/H$ has the shadowing property.
\end{theorem}
\pf
Fix any $E\in\mathfrak{B}_e$ and take $D\in\mathfrak{B}_e$ such that every $D$-pseudo-orbit can be $E$-shadowed by some point in $G$.
Assume that $\{\pi(x_n)\}_{n=1}^{\infty}=\{x_nH\}_{n=1}^{\infty}$ is a $\pi(D)$-pseudo-orbit of $\tilde{f}$. Then,
$\tilde{f}(x_nH)(x_{n+1}H)^{-1}\in \pi(D)$, implying that $f(x_n)x_{n+1}^{-1}H\in DH$. Thus, $f(x_n)x_{n+1}^{-1}\in DH$.
This implies that $\{x_n\}_{n=1}^{\infty}$ is a $DH$-pseudo-orbit of $f$ and by Lemma \ref{lem1} there exists $x\in G$ such that
$f^n(x)x_n^{-1}\in EH$. This implies that $\tilde{f}^n(\pi(x))\pi(x_n)^{-1}\in \pi(E)$.
\epf

\begin{corollary}
Let $f:G\rightarrow G$ be a continuous automorphism with the shadowing property. Then $\tilde{f}:G/H\rightarrow G/H$ has the shadowing property for any choice of H from recurrent subgroups in Proposition \ref{p1}.
\end{corollary}

\begin{proposition}
Let $f:G\rightarrow G$ be a continuous automorphism with the shadowing property. Then, $\Omega(f)=\mathrm{CR}(f)$.
\end{proposition}
\pf
Clearly, $\Omega(f)\subset \mathrm{CR}(f)$. Suppose that $x\in \mathrm{CR}(f)$ and $U$ is an open neighborhood of $x$.
Choose some $E\in\mathfrak{B}_e$ such that $Ex\in U$. By the shadowing property there exists $D\in\mathfrak{B}$ such that
every $D$-pseudo-orbit is $E$-shadowed by some point in $G$. From $x\in \mathrm{CR}(f)$, it follows that there exists a
$D$-pseudo-orbit $\{x_0=x,x_1,\dots,x_n=x\}$ from $x$ to itself. If extend this sequence to an infinite $D$-pseudo-orbit, then this full
pseudo-orbit is $E$-shadowed by some point $z\in G$. Thus, $zx^{-1}, f(z)x_1^{-1},\dots,f^n(z)x^{-1}\in E$, implying that
$z,f^n(z)\in Ex\subset U$. Therefore, $f^n(U)\cap U\neq \O$ and so $x\in \Omega(f)$.
\epf

\section{Topological Entropy}
Let $f$	be a continuous	endomorphism on the	topological group $G$ and
$K$ be a compact subset of $G$. Given $E\in\mathfrak{B}_e$ and $n\in \mathbb{N}$,
a subset $A\subseteq K$ is called an {\it $(n,E,f)$-spanning set} for $K$ if
$$
K\subseteq \bigcup_{x\in A}\left(\bigcap_{i=0}^{n-1}F^{-i}(E)\right)x,
$$
or equivalently, if for any $x\in K$, there exists $y\in A$ such that $f^i(x)f^i(y^{-1})\in E$ for all $i=0,1,2,\dots,n-1$.
By compactness, there exist a finite $(n,E,f)$-spanning set for $K$. Let $\textsf{span}(n,E,K)$ be the minimum cardinality
of all $(n,E,f)$-spanning sets for $K$.

A subset $A\subseteq K$ is called an {\it $(n,E,f)$-separated set} for $K$ if for any pair of distinct points $x$ and $y$ in $A$,
there exists $0\leq i\leq n-1$ such that $f^i(x)f^i(y^{-1})\notin E$. Again by compactness of $K$, every $(n,E,f)$-separated
set for $K$ is finite. Let $\textsf{sep}(n,E,K)$ be the maximum cardinality of all $(n,E,f)$-separated set for $K$.

For any $U\in\mathfrak{B}_e$, define
\begin{eqnarray*}
&&\textsf{span}(E,K)=\limsup_{n\rightarrow\infty}\frac{1}{n}\log \textsf{span}(n,E,K);\\
&&\textsf{sep}(E,K)=\limsup_{n\rightarrow\infty}\frac{1}{n}\log \textsf{sep}(n,E,K).
\end{eqnarray*}
Then, define the following quantities for a uniformly continuous map $f$:
\begin{eqnarray*}
&&h_{\textsf{span}}(f,K)=\sup\left\{\textsf{span}(E,K)~|~E\in\mathfrak{B}_e\right\};\\
&&h_{\textsf{span}}(f)=\sup\left\{h_{\textsf{span}}(f,K)~|~K\in\mathcal{K}(G)\right\};\\
&&h_{\textsf{sep}}(f,K)=\sup\left\{\textsf{sep}(E,K)~|~E\in\mathfrak{B}_e\right\};\\
&&h_{\textsf{sep}}(f)=\sup\left\{h_{\textsf{sep}}(f,K)~|~K\in\mathcal{K}(G)\right\};
\end{eqnarray*}
where $\mathcal{K}(G)$ is the set of all nonempty compact subsets of $G$.

\begin{theorem}\cite{DIKRANJAN20121916}
Let	$f:G\rightarrow G$ be a continuous map on a topological group $G$.
Then, $h_{\textsf{span}}(f)=h_{\textsf{sep}}(f)=h_{top}(f)$.
\end{theorem}

The entropy-carrying sets of a continuous map on a compact space is always in particular interest.
Adapting the techniques in \cite{robinsondynamical}, the following is proved.
\begin{theorem}
Let $f$ be a continuous endomorphism on a Hausdorff compact topological group $G$. Then, $h_{top}(f|_{\Omega(f)})=h_{top}(f)$.
\end{theorem}
\pf
Fix any $m\in \mathbb{N}$ and $E\in\mathfrak{B}_e$, and let $F(m,E,\Omega(f))$ be an $(m,E,f)$-spanning set for $\Omega(f)$ with cardinality $\textsf{span}(m,E,\Omega(f))$. Let $U=\{x\in G~|~\exists y\in F(m,E,\Omega(f)) \text{ such that } f^i(x)f^i(y^{-1})\in E
 \text{ for all } 0\leq i\leq m\}$.

\textbf{Claim 1.} $U$ is an open neighborhood of $\Omega(f)$.

Choose $W\in\mathfrak{B}_e$ such that $W^2\subset E$. Then, there exists $D\in\mathfrak{B}_e$ such that
$xy^{-1}\in D$ implies $f^i(x)f^i(y^{-1})\in W$ for all $0\leq i\leq m$. For any fixed $x\in U$  
and any $y\in Dx$, it is clear that $x^{-1}y\in D$, implying that $f^i(x^{-1})f^i(y)\in W$ for all $0\leq i\leq m$. From $x\in U$,
it follows that there exists $z\in F(m,E,\Omega(f))$ such that $f^i(x)f^i(z^{-1})\in W$ for all $0\leq i\leq m$.
Therefore, $f^i(y)f^i(z^{-1})\in E$. This implies that $Dx\subset U$.

Since $U^c=G\setminus U$ is compact and all points in $U^c$ is wandering, then there exists $V\subset E$ such that
for any $y\in U^c$ and any $n\in \mathbb{N}$, $f^n(Vy)\cap Vy=\O$.
Now let $F(m,V,U^c)$ be an $(m,V,f)$-spanning set for $U^c$ with cardinality $\textsf{span}(m,V,U^c)$ and let
$F_m=F(m,E,\Omega(f))\cup F(m,V,U^c)$. Since $F_m$ is an $(m,E,f)$-spanning set for $G$, we obtain that
$|F_m|\geq \textsf{span}(m,E,G)$. For any $l\in \mathbb{N}$, define
$\phi_l:G\rightarrow F_m^i$ by $\phi_l(x)=(y_0,y_1,\dots, y_{l-1})$, where
\begin{align*}
 & y_j\in F(m,E,\Omega(f))\quad\mbox{  and  }\quad f^{jm}(x)y_j^{-1}\in E\qquad \mbox{ if }\qquad f^{jm}(x)\in U\\
  & y_j\in F(m,V,U^c)\quad\mbox{ and    }\quad f^{jm}(x)y_j^{-1}\in V\qquad \mbox{ if }\qquad f^{jm}(x)\in U^c.
  \end{align*}
  {\color{red}MEANS THAT $ f^{jm}(x)\in U \Rightarrow y_j\in F(m,E,\Omega(f))\cap E^{-1} f^{jm}(x)$}
  
If $\phi_l(x)=(y_0,y_1,\dots, y_{l-1})$ for some $x\in G$, then a point $y_j\in F(m,V,U^c)$ can not be repeated in this $l$-tuple. Because $Ey_j$'s are wandering for any choice of $y_j\in F(m,V,U^c)$.

Choose $n>\textsf{span}(m,V,U^c)$. Let $H(n,E^2,G)$ be an $(n,E^2,f)$-separated set for $G$ with cardinality
$\textsf{sep}(n,E^2,G)$ and let $l$ be a positive integer with $(l-1)m<n\leq lm$.

  \textbf{Claim 2.} The map $\phi_l$ is one to one on $H(n,E^2,G)$.

Suppose that there exists $x,y\in H(n,E^2,G)$ such that $\phi_l(x)=\phi_l(y)=(y_0,y_1,\dots,y_{l-1})$.
For $0\leq i< m$ and $0\leq j <l$, we have
  $f^{i+jm}(x)f^{i+jm}(y^{-1})=f^{i+jm}(x)f^{i}(y_j^{-1})f^{i+jm}(y^{-1})f^{i}(y_j)\in E^2$.
  By the choice of $l$, it follows that for any $0\leq i <n$, $f^{i}(x)f^{i}(y^{-1})\in E^2$.
  This, together with $x,y\in H(n,E^2,G)$, implies that $x=y$.

\textbf{Claim 3.} Let $p=\textsf{span}(m,E,\Omega(f))$ and $q=\textsf{span}(m,V,U^c)$. Then,
$$
|\phi_l(H(n,E^2,G))|\leq (q+1)! l^q p^l.
$$

Let $I_k$ be the set of $l$-tuples in $\phi_l(H(n,E^2,G))$ such that the numbers of components $y_s$
which belongs to $F(m,V,U^c)$ is $k$. Since $y_k\in F(m,V,U^c)$ can not be repeated in $\phi_l(x)$,
then $k\leq q$. For $I_k$, there exist ${{q}\choose{k}}$ ways of picking these $k$ points $y_k\in F(m,V,U^c)$,
there exist $\frac{l!}{(l-k)!}$ ways of arranging these choice among the points in the ordered $l$-tuples.
Meanwhile, there exist at most $p^{l-k}$ ways of picking the remaining $y_s$ from $F(m,E,\Omega)$. Thus,
$$
|I_k|\leq {{q}\choose{k}}\frac{l!}{(l-k)!}p^l.
$$
From ${{q}\choose{k}}\leq q!$ and $\frac{l!}{(l-k)!}\leq l^k$, it follows that
$$
|\phi_l(H(n,E^2,G))|\leq \sum_{k=0}^{q}{{q}\choose{k}}\frac{l!}{(l-k)!}{p^l}\leq (q+1)!l^qp^l.
$$

Now, applying claims 2 and 3 yields that
$$
\textsf{sep}(n,E^2,G)=|\phi_l(H(n,E^2,G))|\leq (q+1)! l^qp^l,
$$
where $p=\textsf{span}(n,E,\Omega(f))$ and $q=\textsf{span}(m,V,U^c)$. This implies that
\begin{eqnarray*}
h_{\textsf{sep}}(f)&\leq &\limsup_{\substack{E\in\mathfrak{B}_e}}\limsup_{\substack{n\rightarrow\infty}}\frac{1}{n}\log \textsf{sep}(n,E^2,G)\\
&\leq& \limsup_{\substack{E\in\mathfrak{B}_e}}\limsup_{\substack{l\rightarrow\infty}}\frac{1}{(l-1)m}[\log ((q+1)!)+q\log(l)+l\log (p)]\\
&\leq& \limsup_{\substack{E\in\mathfrak{B}_e}}\frac{1}{m}\limsup_{\substack{n\rightarrow\infty}}\log p\\
&=&\limsup_{\substack{E\in\mathfrak{B}_e}}\frac{1}{m}\log \textsf{span}(n,E,\Omega).
\end{eqnarray*}
This implies that $h_{top}(f)=h_{\textsf{sep}}(f)\leq h_{top}(f|_{\Omega(f)})$.
\epf
\begin{corollary}
Let $f:G\rightarrow G$ be a continuous endomorphism. Then, $h_{top}(f|_{\mathrm{CR}(f)})=h_{top}(f)$.
\end{corollary}

Addition Theorem states that the entropy is additive in appropriate sense with respect to
invariant subgroups \cite{DIKRANJAN2016612}. Bruno and Virili \cite{bruno2017topological}
proved that the addition theorem in the case of locally compact totally disconnected topological groups.

\begin{theorem}\cite{bruno2017topological}
Let $G$ be a  locally
compact totally disconnected group, and $f:G\rightarrow G$ be a continuous endomorphism, and $H$
be a compact $f$-invariant subgroup of $G$. Then,
$$
h_{top}(f)\geq h_{top}(\tilde{f})+h_{top}(f|_{H}).
$$
\end{theorem}

\begin{corollary}
Let $G$ be a  locally
compact totally disconnected group, and $f:G\rightarrow G$ a continuous endomorphism,
and $\tilde{f}:G/\mathrm{CR}(f)\rightarrow G/\mathrm{CR}(f)$ be the map induced by $f$. Then,
$h_{top}(\tilde{f})=0$.
\end{corollary}

%
%


\section*{Acknowledgment}
Some of the results included in this paper were obtained during  visit of S.A. Ahmadi at the Abdus Salam International Centre for Theoretical Physics (ICTP). Support of this institution is widely acknowledged.

\end{document}